\documentclass[11pt]{amsart}
\usepackage{amssymb, amsmath, amsthm}

\newtheorem{theorem}{Theorem}[section]
\newtheorem*{theorem*}{Theorem}

\newtheorem*{corollary*}{Corollary}
\newtheorem{lemma}[theorem]{Lemma}
\newtheorem{proposition}[theorem]{Proposition}
\newtheorem{definition}[theorem]{Definition}

\newtheorem{problem}[theorem]{Problem}

\newcommand{\N}{\mathbb{N}}

\newcommand{\R}{\mathbb{R}}
\newcommand{\C}{\mathbb{C}}

\newcommand{\ie}{i.e. }

\newcommand{\cstar}{C$^*$}

\newcommand{\continuous}{\mathrm{C}}

\newcommand{\mat}{\mathrm{M}}

\renewcommand{\epsilon}{\varepsilon}
\renewcommand{\phi}{\varphi}
\renewcommand{\emptyset}{\varnothing}
\renewcommand{\leq}{\leqslant}
\renewcommand{\geq}{\geqslant}

\newcommand{\nnegint}{\N}

\newcommand{\posint}{\N^*}

\begin{document}

\title{On perturbations of continuous maps}
\author{Beno\^it Jacob}
\address{University of Toronto\\
Dept. of Mathematics\\
40 St. George Street, room 6290\\
Toronto, Ontario M5S 2E4\\
Canada}
\email{bjacob@math.toronto.edu}

\begin{abstract}
We give sufficient conditions for the following problem: given a topological space $X$, a metric space $Y$, a subspace $Z$ of $Y$, and a continuous map $f$ from $X$ to $Y$, is it possible, by applying to $f$ an arbitrarily small perturbation, to ensure that $f(X)$ does not meet $Z$? We also give a relative variant: if $f(X')$ does not meet $Z$ for a certain subset $X'\subset X$, then we may keep $f$ unchanged on $X'$. We also develop a variant for continuous sections of fibrations, and discuss some applications to matrix perturbation theory.
\end{abstract}

\maketitle

\tableofcontents

\section{Introduction}

Given a topological space $X$, a metric space $Y$, a subspace $Z$ of $Y$, and a continuous map $f$ from $X$ to $Y$, is it possible, by applying to $f$ an arbitrarily small perturbation, to ensure that $f(X)$ does not meet $Z$? And, with respect to which topology on the set of maps from $X$ to $Y$?\\

Sufficient conditions have previously been worked out, with respect to the uniform topology, for cases where $X$ is a CW-complex, using a transversality argument. This then allows, by using projective limit decompositions, to obtain results for all compact Hausdorff $X$. See, for instance, the proof of Theorem 4 in \cite{CE90}, the proof of Lemma 2.5 in \cite{P91}, and the proof of Theorem 3.3 in \cite{P94}. The results of the present article are more general and allow to recover them, see Section \ref{applitomatrix}.\\

This article proceeds from two basic ideas. First, the \emph{source limitation topology} is much more suitable to the study of this problem than the uniform topology is. Second, general topological dimension-theoretic techniques apply here, so there is no need for a reduction to CW-complexes and for transversality arguments.\\

The source limitation topology is the variant of the uniform topology where $\epsilon$ is allowed to vary continuously. Thus we are asking the following question, keeping the notation from above and letting $d$ be the metric on $Y$: given a continuous function $\epsilon:X\rightarrow\R_{>0}$, does there exist a continuous map $g:X\rightarrow Y$ such that for all $x\in X$, $g(x)\not\in Z$ and $d(f(x),g(x))<\epsilon(x)$?\\

In order to simplify subsequent statements, let us introduce the following terminology:
\begin{definition}\label{avoidable}
Let $X$ be a topological space, let $Y$ be a metric space with metric denoted by $d$, and let $Z$ be a subset of $Y$. One says that $Z$ is \emph{$X$-avoidable} if for any continuous map $f:X\rightarrow Y$ and any continuous function $\epsilon:X\rightarrow\R_{>0}$, there exists a continuous map $g:X\rightarrow Y$ such that for all $x\in X$, $g(x)\not\in Z$ and $d(f(x),g(x))<\epsilon(x)$.
\end{definition}

It is worth noting that being $X$-avoidable is a property of the pair $(Z,Y)$, not an intrinsic property of $Z$; nevertheless, when the context doesn't make this ambiguous, we will let $Y$ be implicit and refer to $Z$ alone.\\

Thus our main problem can be rephrased as:
\begin{problem}\label{theproblem}
Given a topological space $X$ and a metric space $Y$, which subsets of $Y$ are $X$-avoidable?
\end{problem}

To fix the terminology, let us take this definition:
\begin{definition}
A \emph{manifold} of dimension $n\in\N$ is a Hausdorff space where every point has a neighborhood that is homeomorphic to $\R^n$.
\end{definition}
Here is our main result:
\begin{theorem*}[See Theorem \ref{avoidsubset}]
Let $n,q\in\nnegint$ with $n<q$. Let $X$ be a normal space of Lebesgue covering dimension $n$. Let $Y$ be a metric topological manifold. Let $Z$ be a subset of $Y$ of \emph{Lipschitz codimension} (see Definition \ref{lipschitzcodim}) at least $q$. Then $Z$ is $X$-avoidable.
\end{theorem*}
Notice that here, while $Y$ is required to be a metric manifold, much greater generality is allowed for $X$.\\

We also give results showing that unions of $X$-avoidable sets are still $X$-avoidable: first a general result on finite unions (Proposition \ref{avoidfiniteunionofsubsets}) and then, in the locally compact case, a result on countably infinite unions (Theorem \ref{avoidunionofsubsets}).\\

We then obtain a relative variant (Proposition \ref{avoidsubsetrelative}):  if the restriction of the map $f$ to a certain closed subset $C$ of $X$ already avoids $Z$, then we may perturb $f$ outside of $C$ to make it avoid $Z$ everywhere, without modifying $f$ on $C$.\\

We also prove a variant for continuous sections of locally trivial fibrations (Theorem \ref{avoidsubfibration}) under the assumption that the base space is paracompact.\\

Finally, in Section \ref{applitomatrix} we discuss some applications to matrix perturbation theory, which motivated this work.

\section{Review of some dimension theory}

Here we review just a few selected notions and results, with no aim to offer a general review of this topic. See \cite{E78} or \cite{N64} for the general theory. Let us first recall the classical notion of Lebesgue covering dimension, which we will just call dimension. Given a topological space $X$, a \emph{refinement} of an open covering $(U_i)_{i\in I}$ of $X$ is an open covering $(V_j)_{j\in J}$ of $X$ such that for all $j\in J$, there exists $i\in I$ such that $V_j\subset U_i$.
\begin{definition}[See \cite{E78} or \cite{N64}]
Let $d\in\nnegint$. A topological space $X$ is said to have \emph{dimension at most $d$} if any open covering of it has a refinement $(V_j)_{j\in J}$ such that for all $x\in X$, the set $\{j\in J,\;x\in V_j\}$ has at most $d+1$ elements.
\end{definition}
Obviously, the dimension of $X$ is then defined as the smallest $d$ such that $X$ has dimension at most $d$, or $\infty$ if no such $d$ exists. It is true that $\R^d$ has dimension $d$.
\begin{definition}
For any topological space $X$, we will let $\dim X$ denote its (Lebesgue covering) dimension.
\end{definition}

Let us also recall the notion of a \emph{normal space}:
\begin{definition}
A topological space is said to be \emph{normal} if any two disjoint closed subsets have disjoint neighborhoods.
\end{definition}

Recall the following classical theorem in dimension theory:
\begin{theorem}[See \cite{N64}, Theorem VII.9]
\label{characdimavoid}
Let $X$ be a normal space. Let $n\in\nnegint$. The following are equivalent:
\begin{enumerate}
\item \label{cond_dimXleqd} $\dim X \leq n$.
\item \label{cond_avoidpointincube} For any $f\in\continuous(X,[0;1]^{n+1})$, for any $\epsilon>0$, for any $y\in[0;1]^{n+1}$, there exists $g\in\continuous(X,[0;1]^{n+1})$ such that $\Vert f-g\Vert <\epsilon$ and $y\not\in g(X)$.
\end{enumerate}
\end{theorem}

We will need some variants and refinements of the implication \ref{cond_dimXleqd}$\Rightarrow$\ref{cond_avoidpointincube} in the above theorem. In order to obtain them, we will simply adapt the classical proof of that theorem. The main technical lemmas used in that proof are the following:

\begin{lemma}[Urysohn's Lemma]\label{urysohnlemma}
Let $X$ be a normal space. Let $F,G$ be disjoint closed subsets of $X$. There exists a continuous function $\phi\in\continuous(X,[0;1])$ such that $\phi(F)=0$ and $\phi(G)=1$.
\end{lemma}
Here, by $\phi(F)=\lambda$ we mean that $\phi(x)=\lambda$ for all $x\in F$.

\begin{lemma}[See \cite{N64}, VII.4.B]\label{bigdimthlemma}
Let $n\in\nnegint$. Let $X$ be a normal space such that $\dim X\leq n$. Let $U_1,\ldots,U_{n+1}$ be open subsets of $X$. Let $F_1,\ldots,F_{n+1}$ be closed subsets of $X$. Suppose that $F_i\subset U_i$ for all $i$. Then there exist open subsets $V_1,\ldots,V_{n+1}$ and $W_1,\ldots,W_{n+1}$ of $X$ such that
$$F_i\subset V_i\subset \overline V_i\subset W_i\subset U_i \;\;\;\text{for}\;i=1,\ldots,n+1$$
and
$$\bigcap_{i=1}^{n+1}\left(\overline W_i - V_i\right)\;=\;\emptyset.$$
\end{lemma}

\section{Perturbations of continuous maps}

The proof of the following lemma follows very closely the classical proof of Theorem \ref{characdimavoid}. Our only change is to replace the uniform topology by the source limitation topology, \ie replace the constant $\epsilon$ by a function on $X$. This is a straightforward adaptation, but will allow us to do much more than the classical theorem allowed.

\begin{lemma}\label{avoidzero}
Let $n\in\nnegint$. Let $X$ be a normal space such that $\dim X\leq n$. For any $f\in\continuous(X,\R^{n+1})$ and for any $\epsilon\in\continuous(X,\R_{>0})$, there exists $g\in\continuous(X,\R^{n+1})$ such that for all $x\in X$, $\Vert g(x)-f(x)\Vert<\epsilon(x)$ and $g(x)\neq 0$.
\end{lemma}
\begin{proof}
Write $f=(f_1,\ldots,f_{n+1})$ where the $f_i\in\continuous(X,\R)$ are continuous functions. For $1\leq i\leq n+1$, let
\begin{eqnarray*}
F_i & = & \{x\in X;\;f_i(x)\geq\epsilon(x)\}\\
G_i & = & \{x\in X;\;f_i(x)\leq-\epsilon(x)\}.
\end{eqnarray*}
Since $F_i\subset X-G_i$ for all $i$, it follows from Lemma \ref{bigdimthlemma} that there exist open subsets $V_1,\ldots,V_{n+1}$ and $W_1,\ldots,W_{n+1}$ of $X$ such that
$$F_i\subset V_i\subset \overline V_i\subset W_i\subset X-G_i \;\;\;\text{for}\;i=1,\ldots,n+1$$
and
\begin{equation}\label{emptycap}\bigcap_{i=1}^{n+1}\left(\overline W_i - V_i\right)\;=\;\emptyset.\end{equation}
Since $X$ is normal, by Urysohn's Lemma \ref{urysohnlemma}, for all $i$ there exists a continuous function $\phi_i:X\rightarrow[-1;1]$ such that $\phi_i(\overline V_i)=1$ and $\phi_i(X-W_i)=-1$. Since $F_i\subset V_i$, we have $\phi_i(F_i)=1$. Since $G_i\subset X-W_i$, we have $\phi_i(G_i)=-1$. We may therefore define a continuous function $g_i:X\rightarrow\R$ by letting, for all $x\in X$,
$$g_i(x)=\left\{\begin{array}{ll}f_i(x) & \text{if}\;x\in F_i\cup G_i \\ \epsilon(x)\phi_i(x) & \text{if}\;x\not\in F_i\cup G_i \end{array}\right.$$
We now define $g:X\rightarrow\R^{n+1}$ by letting $g(x)=(g_1(x),\ldots,g_{n+1}(x))$.
It is clear that $\Vert g_i(x) - f_i(x)\Vert\leq 2\epsilon(x)$ for all $x$ and all $i$, and therefore
$$\Vert g(x)-f(x)\Vert \leq 2\sqrt{n+1}\,\epsilon(x)\;\;\;\text{for all}\;x\in X.$$
It remains to show that $g(x)\neq 0$ for all $x\in X$. Suppose that $g(x)=0$ for some $x\in X$. Then $g_i(x)=0$ for $i=1,\ldots,n+1$. It follows that $x\not\in F_i\cup G_i$, so that $g_i(x)=\epsilon(x)\phi_i(x)=0$. Since $\epsilon(x)>0$, it follows that $\phi_i(x)=0$. This in turn entails that $x\not\in\overline V_i$ and $x\not\in X-W_i$. Therefore, $x\in W_i-\overline V_i$ for all $i=1,\ldots,n+1$, contradicting equation (\ref{emptycap}).
\end{proof}

Let us now introduce the notion of ``Lipschitz codimension''.

\begin{definition}
\label{lipschitzcodim}
Let $p,q\in\nnegint$. Let $Y$ be a metric topological manifold of dimension $p$. Let $Z$ be a subset of $Y$. One says that $Z$ has \emph{Lipschitz codimension at least $q$ in $Y$} if there exists an open covering $(U_i)_{i\in I}$ of $Y$ and, for all $i\in I$, a homeomorphism $\phi_i$ from $U_i$ to an open subset of $\R^p$, such that:
\begin{itemize}
\item for all $i\in I$, $\phi_i^{-1}$ is Lipschitz;
\item for all $i\in I$, $\phi_i(U_i\cap Z)$ is a subset of $\R^{p-q}$ seen as the vector subspace of $\R^p$ consisting of vectors ending with $q$ zeros.
\end{itemize}
\end{definition}

The next step is to prove that we can locally avoid $Z$, that is separately in each open set $U_i$ as in Definition \ref{lipschitzcodim}:

\begin{lemma}
\label{avoidsubsetonechart}
Let $n,p,q\in\nnegint$ with $q>n$. Let $X$ be a normal space of dimension $n$. Let $Y$ be a metric topological manifold of dimension $p$, with
metric denoted by $d$. Let $Z$ be a subset of $Y$ of Lipschitz codimension at least $q$. Let $(U_i)_{i\in I}$ be an open covering of $Y$ and let $(\phi_i)_{i\in I}$ be a family as in Definition \ref{lipschitzcodim}. Let $\epsilon\in\continuous(V_i,\R_{>0})$. Let $f\in\continuous(X,Y)$. Let $i\in I$. Let $V_i=f^{-1}(U_i)$. There exists $g_i\in\continuous(X,Y)$ such that:
\begin{itemize}
\item for all $x\in X-V_i$, $g_i(x)=f(x)$;
\item for all $x\in V_i$, $d(f(x),g_i(x))<\epsilon(x)$;
\item for all $x\in V_i$, $g_i(x)$ does not belong to the closure of $Z$.
\end{itemize}
\end{lemma}
\begin{proof}
Define a function $\eta$ on $V_i$ as follows. If $\phi_i(U_i)\neq\R^p$, let $$\eta(x)=\inf_{a\in\R^p-\phi_i(U_i)} \Vert \phi_i(f(x)) - a\Vert\;\;\text{for all}\;x\in V_i,$$
and if $\phi_i(U_i)=\R^p$, let $\eta(x)=1$ for all $x\in V_i$. Let $\pi:\R^p\rightarrow\R^q$ be the map discarding the $p-q$ first components. Let $\alpha=\pi\circ \phi_i\circ f\vert_{V_i}$. Thus $\alpha$ is a map from $V_i$ to $\R^q$, and we have $\alpha(x)=0$ for all $x\in f^{-1}(Z)\cap V_i$. Let $K_i$ be a Lipschitz constant for $\phi_i^{-1}$.
For $x\in V_i$, let
$$\epsilon'(x)=\min\left(\frac{\epsilon(x)}{K_i},\: \eta(x)\right).$$
By Lemma \ref{avoidzero} applied to $\alpha$ and the function $\epsilon'$, there exists a continuous map $\beta$ from $V_i$ to $\R^q$ such that for all $x\in V_i$, $$\Vert \alpha(x) - \beta(x) \Vert<\epsilon'(x)$$ and $\beta(x)\neq 0$.
Let $\rho:\R^p\rightarrow\R^{p-q}$ be the map discarding the $q$ last components. Define a map $\gamma:V_i\rightarrow\R^p=\R^{p-q}\times\R^q$ by letting $\gamma(x) = (\rho(\phi_i(f(x))), \beta(x))$ for all $x\in V_i$. We have
$$\Vert \gamma(x)-\phi_i(f(x))\Vert = \Vert \beta(x)-\alpha(x)\Vert < \epsilon'(x)$$
for all $x\in V_i$. Since $\epsilon'(x)\leq\eta(x)$, it follows that $\gamma(x)\in\phi_i(U_i)$ for all $x\in V_i$. Thus $\gamma$ is a continuous map from $V_i$ to $\phi_i(U_i)$.\\

Also notice that for all $x\in V_i$, $\gamma(x)$ does not belong to the closure of $\phi_i(Z)$, since the $q$ last components of $\gamma(x)$ are $\beta(x)\neq 0$, and by Definition \ref{lipschitzcodim} we know that any vector in $\phi_i(Z)$ has its last $q$ components equal to 0.\\

Since $\gamma$ is a continuous map from $V_i$ to $\phi_i(U_i)$, we may let $g_i=\phi_i^{-1}\circ \gamma$. Thus $g_i$ is a continuous map from $V_i$ to $U_i$ such that for all $x\in V_i$, $d(f_i(x),g_i(x))<K_i\epsilon'(x)< \epsilon(x)$. Also note that for all $x\in V_i$, since $\gamma(x)$ does not belong to the closure of $\phi_i(Z)$, it follows that $g_i(x)$ does not belong to the closure of $Z$. Finally, by definition of the function $\epsilon'$, we may extend $g_i$ to a continuous map from $X$ to $Y$ by letting $g_i(x)=f(x)$ for all $x\not\in V_i$.
\end{proof}

Let us now glue the local charts together to prove the main theorem of this section:

\begin{theorem}
\label{avoidsubset}
Let $n,q\in\nnegint$ with $q>n$. Let $X$ be a normal space of dimension $n$. Let $Y$ be a metric topological manifold. Let $Z$ be a subset of $Y$ of Lipschitz codimension at least $q$. Then $Z$ is $X$-avoidable.
\end{theorem}
\begin{proof}
Let $d$ be the metric on $Y$. Let $\epsilon\in\continuous(X,\R_{>0})$. Let $f\in\continuous(X,Y)$. We have to show that there exists $g\in\continuous(X,Y)$ such that for all $x\in X$, $d(f(x),g(x))<\epsilon(x)$ and $g(x)$ does not belong to the closure of $Z$. Let $(U_i)_{i\in I}$ be an open covering of $Y$ and $(\phi_i))_{i\in I}$ be a family as given by Definition \ref{lipschitzcodim}, as $Z$ has codimension at least $q$ in $Y$. Since $Y$ is metrizable, it is paracompact, and therefore we may assume without loss of generality that the covering $(U_i)_{i\in I}$ is locally finite --- indeed, the existence of the corresponding family $(\phi_i))_{i\in I}$ in Definition \ref{lipschitzcodim} passes to refinements. Again since $Y$ is paracompact, there is a partition of unity $(u_i)_{i\in I}$ subordinate to $(U_i)_{i\in I}$, and we may replace $U_i$ by the support of $u_i$ so that $(u_i)_{i\in I}$ is precisely subordinate to $(U_i)_{i\in I}$. For $i\in I$, let $V_i=f^{-1}(U_i)$, and let $v_i = u_i\circ f$, so that $(V_i)_{i\in I}$ is a locally finite open covering of $X$ and $(v_i)_{i\in I}$ is a partition of unity precisely subordinate to it. Let $\leq$ be a well-order on $I$. We assume without loss of generality that $(I,\leq)$ is an ordinal. For $x\in X$ and $i\in I$, let
$$\epsilon_i(x)=\sum_{j\leq i} v_i(x)\epsilon(x).$$
Let us prove by transfinite induction that for any ordinal $i$ with $i\leq I$, letting
$$X_i=\bigcup_{j\leq i} V_j,$$
there exists $g_i\in\continuous(X,Y)$ such that $g_i(x)=f(x)$ for all $x\not\in X_i$, and $d(f(x),g_i(x))<\epsilon_i(x)$ and $g_i(x)\not\in Z$ for all $x\in X_i$.\\

\emph{Zero case:} The case $i=0$ follows immediately from Lemma \ref{avoidsubsetonechart}, as we have $X_0=V_0$.\\

\emph{Successor/limit case:} Suppose now that the induction hypothesis is known to hold for all $j<i$ where $i$ is a fixed ordinal, $0<i\leq I$. Let us show that it also holds for $i$. For all $j<i$ we have a map $h_j$ as given by the induction hypothesis. Recall that $(V_j)_{j\in I}$ is a locally finite covering of $X$, so in particular, for all $x\in X$, the set of all $j\in I$ such that $x\in V_j$ is finite. For all $x\in X$, let $j_x\in I$ be the greatest $j<i$ such that $x\in V_j$, and let
$$h_i(x)=g_{j_x}(x).$$
This defines a map $h_i$ from $X$ to $Y$. Let us check that it is continuous.
Again by local finiteness of the covering $(V_j)_{j\in I}$, for any $j<i$, for any $x\in V_j$, there is a neighborhood $W$ of $x$ that intersects only finitely many of the $V_k$, for $k\in I$. Letting $k$ be the greatest among these finitely many indices, we see that $h_i$ agrees with $g_k$ on $W$, hence is continuous at $x$. Thus the map $h_i$ is continuous. Apply Lemma \ref{avoidsubsetonechart} to the map $h_i$, the local chart $(U_i,\phi_i)$, and the epsilon-function $v_i\epsilon$. Call $g_i$ the resulting map. It is then immediate to check that $g_i$ has the wanted properties, showing that our induction hypothesis holds for $i$.
\end{proof}

The following auxiliary results allow to show that certain unions of avoidable subsets are avoidable. They illustrate again how the source limitation topology is more suitable than the uniform topology here.\\

Let us start with finite unions:
\begin{proposition}
\label{avoidfiniteunionofsubsets}
Let $X$ be a normal space. Let $Y$ be a metric topological manifold. Let $n\in\posint$. Let $(Z_i)_{1\leq i\leq n}$ be a family of $X$-avoidable closed subsets of $Y$. Then their union $\bigcup_{1\leq i\leq n} Z_i$ is $X$-avoidable.
\end{proposition}
\begin{proof}
By induction on $n$. The case $n=1$ is vacuously true. Suppose that the results holds for a fixed $n$ and let us establish it for $n+1$. Let $f\in\continuous(X,Y)$ and $\epsilon\in\continuous(X,\R_{>0})$. We may assume that $f$ already avoids $Z_1,\ldots,Z_n$. Let us show that it avoids $Z_{n+1}$. Define a function $\eta$ on $X$ by:
$$\eta_k(x)=\min(\epsilon(x), \min_{i\leq n} d(f(x), Z_i))\;\;\;\text{for all}\;x\in X.$$
Notice that $\eta(x)>0$ for all $x\in X$, because the $Z_i$ are closed. Apply Theorem \ref{avoidsubset} to the map $f$, the subset $Z_{n+1}$ and the function $\eta$. The resulting map $g$ avoids $Z_{n+1}$ by construction, and it still avoids $Z_1,\ldots,Z_n$ because of our particular choice of $\eta$.
\end{proof}

Let us now prove that, when $X$ is locally compact, we may actually avoid the union of countably many closed subsets.
\begin{theorem}
\label{avoidunionofsubsets}
Let $X$ be a normal, locally compact space. Let $Y$ be a metric topological manifold. Let $(Z_i)_{i\in\nnegint}$ be a family of $X$-avoidable closed subsets of $Y$. Then their union $\bigcup_{i\in\nnegint} Z_i$ is $X$-avoidable.
\end{theorem}
\begin{proof}
Let $d$ be the metric on $Y$. For $k\in\nnegint$, let $W_k=\bigcup_{i\leq k}Z_i$. Let $f\in\continuous(X,Y)$ and $\epsilon\in\continuous(X,\R_{>0})$. Let us construct a sequence $(g_k)_{k\in\nnegint}$ of continuous maps from $X$ to $Y$ such that:
\begin{itemize}
\item For all $k\in\nnegint$, for all $x\in X$, $g_k(x)\not\in W_k$.
\item For all $x\in X$, $d(f(x),g_0(x))<\epsilon(x)/2$
\item For all $k\in\nnegint$, for all $x\in X$, $d(g_k(x),g_{k+1}(x))<\eta_k(x)$
where we have put:
$$\eta_k(x)=2^{-k-2}\min(\epsilon(x), \min_{i\leq k} d(g_k(x), Z_i))\;\;\;\text{for all}\;x\in X.$$
\end{itemize}
The existence of $g_0$ follows from the $X$-avoidability of $Z_0$. Let us now suppose $g_0,\ldots,g_k$ to be already constructed for some fixed $k$, and let us construct $g_{k+1}$. Notice that $\eta_k(x)>0$ because the $Z_i$ are closed. Apply Theorem \ref{avoidsubset} to the map $g_k$, the function $\eta_k$, and the subset $Z_{k+1}$. The resulting map $g_{k+1}$ has the wanted properties. This completes the construction of the advertised sequence $(g_k)_{k\in\nnegint}$. It follows from the above-listed properties of that sequence that it converges uniformly on every compact subset of $X$ (since the sequence $(2^{-k-2})_{k\in\nnegint}$ is summable). Let $g$ be its limit. Since $X$ is locally compact, it follows that $g$ is a continuous map from $X$ to $Y$. Notice that for all $x\in X$, we have $d(f(x),g(x))<\epsilon(x)$ (since $1/2+\sum_{k\in\nnegint} 2^{-k-2}=1$). It remains to show that $g$ actually avoids $\bigcup_{i\in\nnegint} Z_i$. Let $x\in X$ and $i\in\nnegint$. We want to show that $g(x)\not\in Z_i$. We have $g_i(x)\not\in Z_i$. For $j\geq i$, let $\delta_j=d(g_j(x),Z_i)$. We have $\delta_j>0$ since $Z_j$ is closed, and by construction we know that $\delta_{j+1}\geq (1-2^{-j-2})\delta_j$ for all $j\geq i$. Since
$$\prod_{j\geq i} 1 - 2^{-j-2} > 0,$$
it follows that $d(g(x),Z_i)>0$, and in particular $g(x)\not\in Z_i$.
\end{proof}

\section{Relative variant}

One of the benefits of the choice of the source limitation topology is that we can easily give a relative version of the above results, \ie if $f$ already avoids $Z$ on a certain closed subset of $X$, we can choose $g$ to agree with $f$ on that subset. However, we will have to make the further assumption that the base space is perfectly normal.\\

Recall that a space $X$ is \emph{perfectly normal} if for any disjoint closed subsets $E,F$ of $X$, there exists a continuous function $\phi$ from $X$ to $[0;1]$ such that $(\phi(x)=0\Leftrightarrow x\in E)$ and $(\phi(x)=1\Leftrightarrow x\in F)$. By contrast, if $X$ is only normal, then the function $\phi$ given by Urysohn's lemma may take the values $0$ or $1$ outside of $E$ and $F$. For instance, metrizable spaces are perfectly normal.

\begin{proposition}\label{avoidsubsetrelative}
Let $X$ be a perfectly normal space. Let $C$ be a closed subset of $X$. Let $X'$ be the complement of $C$ in $X$. Let $Y$ be a metric topological manifold, with metric $d$. Let $Z$ be a $X'$-avoidable subset of $Y$. Let $f$ be a continuous map from $X$ to $Y$. Let $\epsilon$ be a continuous function from $X$ to $\R_{>0}$. Suppose that
$$\text{for all}\;\;x\in C, \;\;\;f(x)\not\in Z.$$
Then there exists a continuous map $g$ from $X$ to $Y$ such that for all $x\in C$, $g(x)=f(x)$ and for all $x\in X$, $g(x)\not\in Z$ and $d(f(x),g(x))<\epsilon(x)$.
\end{proposition}
\begin{proof}
Since $X$ is perfectly normal and $C$ is closed, there exists a continuous function $\eta$ from $X$ to $[0;1]$ vanishing exactly on $C$. Let $\epsilon'$ be the restriction of $\eta\epsilon$ to $X'$, and let $f'$ be the restriction of $f$ to $X'$. Since $Z$ is $X'$-avoidable, there exists a continuous map $g'$ from $X'$ to $Y$ such that for all $x\in X'$, $g'(x)\not\in Z$ and $d(f'(x),g'(x))<\epsilon'(x)$. Finally, extend $g'$ into a continuous map $g$ from $X$ to $Y$ by letting $g(x)=f(x)$ for all $x\not\in X'$.
\end{proof}

\section{Variant for locally trivial fibrations}\label{fibrations}

In this section, we extend our results to the case of locally trivial fibration, under the additional assumption that the base space $X$ is paracompact. Recall that paracompact spaces are normal (Dieudonn\'e's theorem). The proof is similar to that of Theorem \ref{avoidsubset}, except that we now work on the base space $X$ of the fibration instead of working on $Y$ (which is now the fibre of the fibration), whence the need to make the paracompactness assumption on $X$.

\begin{theorem}\label{avoidsubfibration}
Let $X$ be a paracompact space. Let $Y$ be a metric topological manifold, with metric $d$. Let $Z$ be a subset of $Y$ that is $U$-avoidable for all open subsets $U$ of $X$. Let $A$ be a locally trivial fibration over $X$ with fibre $Y$. Let $B$ be a locally trivial sub-fibration of $A$ over $X$ with fibre $Z$. Let $f$ be a continuous section of $A$ over $X$. Let $\epsilon$ be a continuous function from $X$ to $\R_{>0}$. Then:
\begin{itemize}
\item There exists a continuous section $g$ of $A$ over $X$ such that for all $x\in X$, $g(x)\not\in B_x$ and $d(f(x),g(x))<\epsilon(x)$.
\item If moreover $X$ is perfectly normal and $C$ is a closed subset such that for all $x\in C$, $f(x)\not\in B_x$, then the map $g$ may be taken so that for all $x\in C$, $g(x)=f(x)$.
\end{itemize}
\end{theorem}
\begin{proof}
Once the first statement is proved, the proof of the second statement is similar to that of Proposition \ref{avoidsubsetrelative}, so let us only prove the first statement here. Let $\epsilon\in\continuous(X,\R_{>0})$. Let $f$ be a continuous section of $A$ over $X$. We have to show that there exists a continuous section $g$ of $A$ over $X$ such that for all $x\in X$, $d(f(x),g(x))<\epsilon(x)$ and $g(x)$ does not belong to the $B_x$. Since $A$ and $B$ are locally trivial and $Z$ has codimension (Definition \ref{lipschitzcodim}) at least $q$ in $Y$, there exists a covering $(U_i)_{i\in I}$ of $X$ trivializing $A$ and $B$, and there exists a family $(\phi_i))_{i\in I}$ as in Definition \ref{lipschitzcodim} applied to $Z$. Since $X$ is paracompact, we may assume without loss of generality that the covering $(U_i)_{i\in I}$ is locally finite --- indeed, the existence of the corresponding family $(\phi_i))_{i\in I}$ in Definition \ref{lipschitzcodim} passes to refinements. Again since $X$ is paracompact, there is a partition of unity $(u_i)_{i\in I}$ subordinate to $(U_i)_{i\in I}$, and we may replace $U_i$ by the support of $u_i$ so that $(u_i)_{i\in I}$ is precisely subordinate to $(U_i)_{i\in I}$. Let $\leq$ be a well-order on $I$. We assume without loss of generality that $(I,\leq)$ is an ordinal. For $x\in X$ and $i\in I$, let
$$\epsilon_i(x)=\sum_{j\leq i} u_i(x)\epsilon(x).$$
Let us prove by transfinite induction that for any ordinal $i$ with $i\leq I$, letting
$$X_i=\bigcup_{j\leq i} U_j,$$
there exists $g_i\in\continuous(X,Y)$ such that $g_i(x)=f(x)$ for all $x\not\in X_i$, and $d(f(x),g_i(x))<\epsilon_i(x)$ and $g_i(x)\not\in B_x$ for all $x\in X_i$.\\

\emph{Zero case:} The case $i=0$ follows immediately from Theorem \ref{avoidsubset}, as we have $X_0=V_0$ on which the fibrations $A$ and $B$ are trivial.\\

\emph{Successor/limit case:} Suppose now that the induction hypothesis is known to hold for all $j<i$ where $i$ is a fixed ordinal, $0<i\leq I$. Let us show that it also holds for $i$. For all $j<i$ we have a map $h_j$ as given by the induction hypothesis. Recall that $(U_j)_{j\in I}$ is a locally finite covering of $X$, so in particular, for all $x\in X$, the set of all $j\in I$ such that $x\in U_j$ is finite. For all $x\in X$, let $j_x\in I$ be the greatest $j<i$ such that $x\in U_j$, and let
$$h_i(x)=g_{j_x}(x).$$
This defines a map $h_i$ from $X$ to $Y$. Let us check that it is continuous.
Again by local finiteness of the covering $(U_j)_{j\in I}$, for any $j<i$, for any $x\in U_j$, there is a neighborhood $W$ of $x$ that intersects only finitely many of the $U_k$, for $k\in I$. Letting $k$ be the greatest among these finitely many indices, we see that $h_i$ agrees with $g_k$ on $W$, hence is continuous at $x$. Thus the map $h_i$ is continuous. Apply Theorem \ref{avoidsubset} to the map $h_i$ defined on the normal space $U_i$ and the epsilon-function $u_i\epsilon$. Call $g_i$ the resulting map. It is then immediate to check that $g_i$ has the wanted properties, showing that our induction hypothesis holds for $i$.
\end{proof}

\section{Applications to matrix perturbation theory}\label{applitomatrix}

A common pattern of questions is whether certain matrix fields over a certain space may be perturbed to satisfy pointwise a certain condition. For example, may a unitary matrix field $f$ over a space $X$ be perturbed so that at every point $x\in X$, the matrix $f(x)$ has no repeated eigenvalues?\\

This is a special case of Problem \ref{theproblem}. Let us treat an example to illustrate this.\\

Recall these known results by N. Christopher Phillips:
\begin{lemma}[See {\cite[Lemma 2.4]{P91}}]\label{phil1}
The set of elements of $SU(n)$ with at least one repeated eigenvalue is the union of finitely many submanifolds of $SU(n)$, all of codimension at least 3.
\end{lemma}
\begin{lemma}[See {\cite[Lemma 2.5]{P91}}]\label{phil2}
Let $X$ be a finite simplicial complex of dimension at most 2. Let $E$ be a locally trivial $\mat_n(\C)$-bundle over $X$, let $u\in\Gamma(SU_E)$, and let $\epsilon>0$. Then there exists $v\in\Gamma(SU_E)$ such that $\Vert u-v\Vert<\epsilon$ and $v(x)$ has no repeated eigenvalues for all $x\in X$.
\end{lemma}

It follows immediately from Lemma \ref{phil1} and our results that:
\begin{lemma}
Let $Y=SU(n)$. Let $Z$ be the subset of $Y$ of elements with at least one repeated eigenvalue. Then $Z$ is $X$-avoidable for all normal spaces $X$ of dimension at most $2$.
\end{lemma}
\begin{proof}
$Z$ is a finite union of submanifolds of codimension at least $3$, each of which is $X$-avoidable by Theorem \ref{avoidsubset}, so $Z$ is $X$-avoidable by Proposition \ref{avoidfiniteunionofsubsets}.
\end{proof}

Combining this with our results on fibrations, we get the following generalization Lemma \ref{phil2}:
\begin{lemma}\label{generalizephil2}
In Lemma \ref{phil2}, the assumption that $X$ is a finite simplicial complex may be weakened to just the assumption that $X$ any paracompact space. If the $\mat_n(\C)$-bundle $E$ is trivial, then it can be further weakened to just the assumption that $X$ is any normal normal space. Moreover, if $X$ is perfectly normal, and if $u$ is already in the wanted form on a certain closed subset, then $v$ may be taken to agree with $u$ on that subset.
\end{lemma}
\begin{proof}
This is just applying Theorem \ref{avoidsubfibration}.
\end{proof}

\def\cprime{$'$}

\end{document}